\documentclass{article}

\usepackage{amsmath,amssymb,amsthm,color,graphicx}

\newcommand{\gb}{Gr\"obner basis }
\newcommand{\gbs}{Gr\"obner bases }
\def \cocoa{{\hbox{\rm C\kern-.13em o\kern-.07em C\kern-.13em o\kern-.15em A}} }

\newtheorem{thm}{Theorem}[section]

\newtheorem{lem}[thm]{Lemma}

\newtheorem{defn}[thm]{Definition}
\newtheorem{rem}[thm]{Remark}

\begin{document}

\title{Efficiently Computing Gr\"obner Bases of Ideals of Points}

\author{Winfried Just$^{1}$ and Brandilyn Stigler$^2$\\ \\
{\small $^1$Department of Mathematics, Ohio University, Athens, OH 45701}\\
{\small $^2$Mathematical Biosciences Institute, The Ohio State
University, Columbus, OH 43210}}
%
%
\maketitle

\begin{center}
\textsc{\small{Dedicated to Avner Friedman on the occasion
of his 75th birthday}}
\end{center}

\hfill 

\abstract{We present an algorithm for computing \gbs of vanishing
ideals of points that is optimized for the case when the number of
points in the associated variety  is less than the number of
indeterminates. The algorithm first identifies a set of essential
variables, which reduces the time complexity with respect to the
number of indeterminates, and then uses PLU decompositions to reduce
the time complexity with respect to the number of points. This gives
a theoretical  upper bound for its time complexity that is an order
of magnitude lower than the known one for the standard
Buchberger-M\"oller algorithm if the number of indeterminates is
much larger than the number of points. Comparison of implementations
of our algorithm and the standard Buchberger-M\"oller algorithm in
\emph{Macaulay 2} confirm the theoretically predicted speedup. This
work is motivated by recent applications of \gbs to the problem of
network reconstruction in molecular biology.}


\section{Introduction}

Recently, \gbs have been proposed as a promising selection tool
in applications to molecular biology \cite{LS, DJLS}. In these
applications, the data consists of $m$ vectors of discretized
concentration values in a finite field~$k$ for a network of~$n$
biochemicals.  The data points can be viewed as an affine
variety $V$ with points in $k^n$ of multiplicity one and
correspond to the vanishing ideal $\mathbf{I}(V)$ of these
points in the polynomial ring $k[x_1, \ldots , x_n]$.  Each
variable $x_i$ represents the $i$-th biochemical which takes on
values in~$k$.
 Typically, the number of data points $m = |V|$ is on the
order of tens, while the number of variables~$n$ may be in the
thousands (for example, see \cite{yeung}).  This requires
finding \gbs in situations were $m \ll n$, and the run-time of
algorithms for this step constitutes a bottleneck for overall
feasibility of these calculations.  The primary motivation of
this paper is to find an algorithm that optimizes run-time in
the case when $m \ll n$.

Several methods have been described and implemented for computing
\gbs and the associated standard monomials of vanishing ideals of
points. In \cite{BM}, the authors presented the Buchberger-M\"oller
(BM) algorithm for computing the reduced \gb of the ideal of a
variety $V$ over a field. The BM algorithm performs Gaussian
elimination on a generalized Vandermonde matrix and its complexity
is quadratic in the number of indeterminates and cubic in the number
of points in $V$ \cite{abbott, marinari, mora, robbiano}.
 Farr and Gao presented an algorithm based on a
generalization of Newton interpolation \cite{farr}. While the
complexity of their algorithm is exponential in the number $n$ of
indeterminates, the algorithm has been designed for the case in
which $n$ is small as compared to the number of points.
 Lederer proposed a method for lexicographic term
orders which gives insight into the structure of the Gr\"obner basis
\cite{lederer}.  Cerlienco and Mureddu proposed a combinatorial
method that uses Ferrers diagrams to compute the set of standard
monomials for the vanishing ideal of a given set of points with
respect to an inverse lexicographical order \cite{cerlienco}.

In \cite{JS}, the present authors introduced a modification of
BM specifically for the case when the number of points~$m$ in a
given variety is less than the number of indeterminates~$n$.
The EssBM (for Essential Buchberger-M\"oller) algorithm
proposed in that paper identifies \textit{essential} variables,
that is, those in the support of the standard monomials
associated to the ideal of the points, and computes the
relations in the reduced \gb in terms of these variables using
BM.  Since the standard monomials are in terms of at most~$m$
variables, the computation of a \gb can be restricted to a
proper subring of the underlying ring involving only the
essential variables. EssBM was shown to have a worst-case
complexity of $O(nm^3+m^6)$, which is dominated by the first
term when $n\gg m$.

Here we present an improvement of the EssBM algorithm in which we
eliminate the use of BM altogether.  This new algorithm, which we
call EssGB (for Essential Gr\"obner Bases), makes use of PLU
decompositions providing an overall improvement in worst-case
complexity to $O(nm^2 + m^4)$ for a fixed finite field.

The remainder of our paper is organized as follows. In Section 2 we
give a description of the algorithm and in Section 3 we provide the
theoretical background for it. In Section~4 we estimate the
worst-case time complexity of our algorithm.   We conclude with a
summary of the performance of an implementation of our algorithm in
the computer algebra system \emph{Macaulay 2} on test data. These
empirical tests confirm the theoretically predicted speedup relative
to implementations of BM and EssBM on the same platform.

\section{The EssGB Algorithm}

Throughout this paper, let $R=k[x_1,\ldots ,x_n]$ denote a
polynomial ring over a finite field $k$, and let $\prec$ be a fixed
term order on $R$.  For $a=(a_1,\ldots,a_n)\in \mathbb
Z_{\geq 0}^n$, let $x^a$ denote the monomial $x_1^{a_1}\cdots
x_n^{a_n}$.

\begin{defn}
    The \emph{support} of a monomial $x^a \in R$ is $supp(x^a)=\{x_i : x_i | x^a\}$.
\end{defn}

This is not to be confused with the support of a
\emph{polynomial}~$f$, denoted $Supp(f)$, which is the set of
monomials that occur in~$f$.

Let $V\subset k^n$ be a variety of points with multiplicity one
and $|V|=m<\infty$. We consider the problem of computing the
reduced \gb of the vanishing ideal $\mathbf{I}(V)$ of the
points in $V$ with respect to~$\prec$. We call a \gb $G$
\emph{reduced} if all generators are monic (leading
coefficients are equal to 1) and for all $g,h\in G$, if $g\neq
h$, then the leading term of $g$ does not divide any monomial
in $Supp(h)$. We call a polynomial $f\in R$ \emph{reduced with
respect to $G$} if $f$ is the normal form of a polynomial
$f'\in R$ with respect to~$G$, that is,~$f$ is the remainder
of~$f'$ upon division by the elements of~$G$. If the context is
clear, we simply say that $f$ is reduced.

Let $I\subset R$ be an ideal and $G$ a \gb for $I$ with respect to
$\prec$. For any $f\in I$, let $LT(f)$ denote the leading term of
$f$ with respect to $\prec$ and $tail(f)$ the polynomial $f-LT(f)$.
The ideal generated by the set $\{LT(g) : g\in G\}$ is denoted by
$LT(G)$. Further, let $SM(G)$ be the set of monomials not in
$LT(G)$. Note that $SM(G)$ is a $k$-vector space basis for $R/I$.
 We call $SM(G)$ the \emph{set of standard monomials
associated to~$G$}.  In this paper, we restrict our attention to the
case where $I=\mathbf I(V)$ for a finite variety, that is,~$I$ is a
zero-dimensional radical ideal, and $SM(G)$ is a finite basis for
$R/I$.

\begin{defn}
    A variable $x_i$ is \emph{essential} if $x_i\in
    SM(G)$.
\end{defn}

Equivalently, $x_i$ is essential if and only if there is
a monomial $x^a\in SM(G)$ such that $x_i\in supp(x^a)$.
Let $EV(G)$ denote the union of the supports of the standard
monomials $x^a \in SM(G)$. Note that $LT(G)$, $SM(G)$, and $EV(G)$
depend only on the ideal $\mathbf{I}(V)$ and the term order~$\prec$.
Thus we can indicate this dependence by the notation chosen here.

Let $P=\{p_1,\ldots ,p_s\}\subset k^n$ be a set of points. A
polynomial $f\in R$ is a \emph{separator of $p_i\in P$} if
$f(p_i)=1$ and $f(p_j)=0$ for all other $p_j\in P$. Given a variety
$V$ of points of multiplicity one and a term order $\prec$, the
EssGB algorithm returns the triple $(G, SM(G), S)$, where $G$ is the
reduced \gb of the ideal $\mathbf{I}(V)$ of points in $V$ with
respect to~$\prec$; $SM(G)$ is the set of standard monomials
associated to $G$; and $S$ is the set of
reduced separators of the
points in $V$.

Initialize each set as follows: $EV_0=\{\}$ and $SM_0=\{1_R\}$. Let
$[n]$ denote the set $\{1,\ldots ,n\}$ and  for $i\in [n]$, let
$EV_i$ and $SM_i$ denote $i$-th approximations of the corresponding
sets.

For each $i\in [n]$, do the following. Find the $i$-th smallest
variable, say $x_i$. Suppose there are $r$ monomials
$x^{a_1},\ldots,x^{a_r}$ in $SM_{i-1}$. Note that these are
$k$-linearly independent. Try to write $x_i$ as a $k$-linear
combination of these monomials. That is, find (if they exist)
$c_1,\ldots ,c_r\in k$, where
\begin{equation}
\label{linsys}
\begin{split}
  x_i(p_1) &= \sum_{j=1}^r c_jx^{a_j}(p_1) \\
  x_i(p_2) &= \sum_{j=1}^r c_jx^{a_j}(p_2) \\
  &\cdots  \\
  x_i(p_m) &= \sum_{j=1}^r c_jx^{a_j}(p_m)
\end{split}
\end{equation}
and $x^a(p_t)$ is the evaluation of $x^a$ at the $t$-th point in $V$
for $t\in[m]$.

For solving the system~(1) we will use a PLU decomposition
$P_{i-1}L_{i-1}U_{i-1}$ of the matrix $A_{i-1}=(x^{a_j}(p_t))$ of
the monomials $x^{a_j}\in SM_{i-1}$ evaluated at the points in $V$.
This will reduce the time complexity at each step at which no new
essential variable is added to $EV_{i-1}$.  Note that, in general,
$A_{i-1}$ will not be square, but will have dimensions $m \times r$,
where $r = |SM_{i-1}| \leq m$. Still, the standard Gaussian
elimination procedure can be applied to find matrices $P_{i-1},
L_{i-1}, U_{i-1}$ whose product is $A_{i-1}$ and such that $P_{i-1}$
has dimensions $m \times m$ and undoes all row exchanges of the
Gaussian elimination, $L_{i-1}$ is an $m \times m$ lower triangular
matrix with ones on the main diagonal, and $U_{i-1}$ has dimensions
$m \times r$ and is upper triangular in the sense that $u_{k\ell} =
0$ whenever $k > \ell$. Thus even if $A_{i-1}$ is not square, it has
a PLU decomposition in the above sense.

If the system (1) has a solution, then it must be unique (see Lemma
\ref{PLUlemma}). If $c_j = 0$ whenever $x_i \prec x^{a_j}$,
then~$x_i$ is inessential and is the leading monomial of a
polynomial in~$\mathbf I(V)$. In this case, $EV_i = EV_{i-1}$,
$SM_{i} = SM_{i-1}$, $A_i = A_{i-1}$, $P_i = P_{i-1}$, $L_i =
L_{i-1}$, and $U_i = U_{i-1}$.

If no solution exists or $c_j \neq 0$ for some $j$ with $x_i
\prec x^{a_j}$, then $x_i$ is an essential variable and hence
is a standard monomial. In this case let $EV_i=EV_{i-1}\cup
\{x_i\}$; compute the set $SM_{i}$ of standard monomials for
the ideal $\mathbf{I}(V)\cap k[EV_{i}]$ of the points projected
onto the variables in $EV_{i}$ (see Lemma \ref{SMlemma}); and
compute the PLU decomposition $P_{i}L_{i}U_{i} = A_{i}$ of the
matrix $A_i=(x^{a_j}(p_t))$ of the monomials $x^{a_j}\in SM_i$
evaluated at the points in~$V$.

At the end of the loop, all essential variables and standard
monomials have been identified.  The minimum set (with respect to
inclusion) of generators $x^a$ of the leading term ideal of
$\mathbf{I}(V)$ is identified (see Lemma \ref{LTlemma}), and for
each of these generators a polynomial $x^a-g$ is computed so that
the set of all of these polynomials forms a reduced Gr\"obner basis
for $\mathbf{I}(V)$. Finally, the set $S$ of
reduced separators is then computed by solving a system of
linear equations.

\subsection{EssGB}

Let $M$ be an $(m\times n)$-matrix with rows being the points
of $V$, and $\prec$ a term order.  We will assume that
$x_1\prec \cdots \prec x_n$.

\bigskip

\noindent \textbf{Input}: $M$; $\prec$
\smallskip

\noindent \textbf{Output}: $(GB, SM_n,S)$ where $GB$ is the
(reduced) \gb for $\mathbf{I}(V)$ with respect to $\prec$, $SM_n$ is
the set of standard monomials for $GB$, and $S$ is the set of
reduced separators of the points in $V$.
\smallskip

\begin{enumerate}
    \item Initialize $EV_0:=\{\}, SM_0:=\{1_R\}, GB:= \{\}, A_0:=[1,\ldots,1]^T\in k^m,
    P_0 :=Id(m)$, $U_0:=[1,0,\ldots,0]^T\in k^m$, and let $L_0$ be
    the $(m\times m)$-matrix that has ones in the first column and on the diagonal, and zeros elsewhere.
    \item FOR $i\in \{1..n\}$ do
    \begin{enumerate}
        \item Initialize $x_i := i$-th smallest variable, $r:=|SM_{i-1}|$, and $b_i := i$-th column of~$M$.
        \item IF there is no solution
            $c=[c_1,\ldots,c_r]^T$ to the system
            $P_{i-1}L_{i-1}U_{i-1} \cdot c=b_i$ such
            that $c_j=0$ whenever $x_i \prec x^{a_j}$\\
            THEN
        \begin{enumerate}
            \item $EV_{i}:=EV_{i-1}\cup \{x_i\}$.
            \item Compute $SM_i$ in $k[EV_{i}]$ using the algorithm SM-A.
            \item Compute the matrix $A_i:=(x^{a_j}(p_t))$,
            for $x^{a_j}\in SM_i$ and $p_t$ the point in row $t$ of
            $M$.
            \item Compute the PLU decomposition $P_iL_iU_i$ of $A_i$.
        \end{enumerate}
    \end{enumerate}
    \item Compute the set $LT$ of generators of the leading term ideal of $\mathbf I(V)$
    using the algorithm LT-A.
    \item FOR $j\in \{1..|LT|\}$ do
    \begin{enumerate}
        \item Let $b_j=(x^{d_j}(p_t))$ be the $(m\times 1)$-vector of values of
        the monomial $x^{d_j}\in LT$ evaluated at the points $p_t$ in $M$.
        \item Find a solution $[c_1,\ldots,c_m]^T$ of $P_nL_nU_n \cdot c=b_j$.
        \item $GB=GB\cup \{x^{d_j}-\sum c_\ell x^{a_\ell}\}$ where
        $x^{a_\ell}\in SM_n$.
    \end{enumerate}
    \item Compute the set $S$ of reduced separators for $M$ using the algorithm SP-A.
    \item RETURN $GB$, $SM_n$, and $S$.
\end{enumerate}

\subsection{Supporting algorithms}

This section contains the subroutines used in the main algorithm
EssGB.

\subsubsection{SM-A}

The algorithm SM-A generates a set $SM_i$ of standard monomials for
$\mathbf I(V)\cap k[EV_i]$, given a newly identified essential
variable $x_i$ and the set $SM_{i-1}$ of standard monomials for
$\mathbf I(V)\cap k[EV_{i-1}]$. It first constructs a sorted set of
candidate monomials by forming all products $x_i^qx^a$ of monomials
in $SM_{i-1}$ and powers of $x_i$, for $0\leq q < |k|$. Then the
monomials which are $k$-linearly independent can be found by
identifying the pivots of the evaluation matrix
$A_i:=(x^{a_j}(p_t))$, where $x^{a_j}\in C$ and $p_t$ is the point
in row $t$ of $M$.

\bigskip

\noindent \textbf{Input}: $x_i$ an essential variable; $SM_{i-1}$.
\smallskip

\noindent \textbf{Output}: $SM_i$ the set of standard
monomials for $\mathbf I(V)\cap k[EV_i]$.
\smallskip

\begin{enumerate}
    \item Compute the set $C=\{x_i^qx^a:x^a\in
        SM_{i-1},0\leq q < |k|\}$ of candidate standard
        monomials. \item Sort $C$ so that
        $C=\{x^{a_1},\ldots,x^{a_s}:s=|C|, x^{a_j}\prec
        x^{a_{j+1}}\text{ for all }j\}$.
    \item Compute the matrix $A:=(x^{a_j}(p_t))$,  for $x^{a_j}\in C$ and $p_t$ the point in row $t$ of
    $M$.
    \item Compute the row-echelon form $U$ of $A$.
    \item Identify the columns $\pi(1),\ldots,\pi(r)$ corresponding to the $r\leq s$ pivots of~$U$.
    \item RETURN $SM_i=\{x^{a_{\pi(1)}},\ldots,x^{a_{\pi(r)}}\}$.
\end{enumerate}

\subsubsection{LT-A}

This algorithm identifies all minimal leading terms $x^a$ of
$\mathbf I(V)$. We use the following observation, which will be
proved in the next section (Lemma~\ref{LTlemma}).

\begin{rem}
The ideal $LT(G)$ is generated by variables $x_i\notin EV(G)$ and
monomials $x^a$ such that $supp(x^a)\subset EV(G)$, $x^a \notin
SM(G)$, and $x^a$ is  \textit{minimal} in the sense that no monomial
in $LT(G)$ divides $x^a$.
\end{rem}

Recall that $SM_n$ and $EV_n$ are the sets of standard
monomials and essential variables, respectively, after the execution
of Step 2. We will assume that $SM_n$ and $EV_n$ are sorted
according to $\prec$.

\bigskip

\noindent \textbf{Input}: $SM_n$;
$EV_n$.
\smallskip

\noindent \textbf{Output}: $LT$ the set of generators of
the leading term ideal of $\mathbf I(V)$.
\smallskip

\begin{enumerate}
    \item Initialize $C:=$ the $(r\times m)$-matrix of ones, where $r:=|EV_n|, m:=|SM_n|$; $LT:=\{\}$.
    \item FOR $i\in \{1..r\}$ do
    \begin{enumerate}
        \item FOR $j\in \{1..m\}$ do
        \begin{enumerate}
            \item IF $x_ix^{a_j}\in SM_n$ where $x_i\in EV_n$ and
            $x^{a_j}\in SM_n$
            \item THEN $C(i,j):=0$
            \item ELSE FOR $k\in \{1..j-1\}$ do
            \begin{enumerate}
                \item IF $C(i,k)==1$ AND $(x^{a_j})\%(x^{a_k})==0$
                \item THEN $C(i,j):=0$.
            \end{enumerate}
        \end{enumerate}
    \end{enumerate}
    \item FOR $i\in \{1..r\}$ do
    \begin{enumerate}
        \item FOR $j\in \{1..m\}$ do
        \begin{enumerate}
            \item IF $C(i,j)==1$
            \item THEN $LT=LT\cup \{x_ix^{a_j}\}$.
        \end{enumerate}
    \end{enumerate}
    \item Remove repeated elements in $LT$.
    \item $LT = LT\cup \{x_i\not \in EV_n\}$.
    \item RETURN $LT$.
\end{enumerate}

\subsubsection{SP-A}

The algorithm SP-A computes the separators of the points in
$V=\{p_1,\ldots,p_m\}$ in terms of the standard monomials associated
to the ideal of the points. For each point~$p_t$, we wish to find a
polynomial $s_t(\mathbf{x})=\sum_{j=1}^m c_jx^{a_j}\in k[x_1,\ldots
,x_n]$ that satisfies the following:
\begin{eqnarray*}
  \sum_{j=1}^m c_jx^{a_j}(p_t) = 1; & \hspace{0.25in} &
  \sum_{j=1}^m c_jx^{a_j}(p_\ell) = 0,   \text{\hspace{0.1in} for all }\ell \neq t.\\
\end{eqnarray*}
We can do so by solving the system $A_nc=e_t$, where $A_n$ is the
evaluation matrix $A_n=(x^{a_j}(p_\ell))_{\ell,j\in\{1..m\}}$
constructed during execution of EssGB, $c=[c_1,\ldots,c_m]^T\in k^m$
is a vector of unknowns, and $e_t$ is a standard column basis
vector.

\bigskip

\noindent \textbf{Input}: $SM_n=\{x^{a_1}=1,\ldots,x^{a_m}\}$, the
set of standard monomials in increasing $\prec$-order; $(m\times
m)$-matrix $A_n$ in its PLU form
$P_nL_nU_n$.
\smallskip

\noindent \textbf{Output}:
$S=\{s_1(\mathbf{x}),\ldots,s_m(\mathbf{x})\}$ the set of
reduced separators of the points in $V$.
\smallskip

\begin{enumerate}
    \item Initialize $S=\{\}$.
    \item FOR $t\in \{1..m\}$ do
    \begin{enumerate}
        \item Compute $c=[c_1,\ldots,c_m]^T$ such that $P_nL_nU_n\cdot
        c=e_t$.
        \item $S=S\cup \{s_t(\mathbf{x}):=\sum_{j=1}^m c_jx^{a_j}\}$.
    \end{enumerate}
    \item RETURN $S$.
\end{enumerate}

\section{Theoretical Background}

\subsection{SM-A}

Recall that the SM-A algorithm computes a sorted list
$C=\{x^{a_1},\ldots,x^{a_s}\}$ of candidate monomials and returns
$SM_i=\{x^{a_{\pi(1)}},\ldots, x^{a_{\pi(r)}} \}\subset C$, where
$\pi(1),\ldots,\pi(r)$ refer to the columns of the row echelon form
of $A$ corresponding to pivots and $A:=(x^{a_j}(p_t))$ is the
evaluation matrix computed in Step 3 of the subroutine.

\begin{lem}
\label{SMlemma}
    Let $SM_i$ be the output returned by the SM-A subroutine, given
    an essential variable $x_i$ and the set $SM_{i-1}$ of standard
    monomials for $\mathbf I(V)\cap k[EV_{i-1}]$.  Then $SM_i$ is
    the set of standard monomials for $\mathbf I(V)\cap k[EV_i]$, where $EV_i =
    EV_{i-1}\cup \{x_i\}$.
\end{lem}

\begin{proof}
    The set $C$ consists of all multiples of $x_i$ and $x^a\in SM_{i-1}$
    and so generates the $k$-vector space $R/I\cap k[EV_i]$.
    The dimension of this space is equal to the number $r$ of nonzero rows of
    the matrix $U$ as computed in Step~(4) of SM-A.

    Now consider $x^{a_k} \in C$ that is not in the set $SM_i$ returned by SM-A, and let $U(k)$ consist of the
    first $k$ columns of $U$.  Then $U(k-1)$ and $U(k)$ have the same rank and it follows that
    the $k$-th column of $A$ is a linear combination of the columns of $A$ indexed $j = 1, \ldots , k-1$.
    This means that
    \begin{equation}\label{x^k}
    x^{a_{k}} - \sum_{j = 1}^{k-1} c_j x^{{a_j}} \in \mathbf{I}(V),
    \end{equation}
    for some coefficients $c_j$.  Since the elements of $C$ were listed in increasing order
    with respect to $\prec$, the monomial~$x^{a_k}$ is the leading monomial in~(\ref{x^k}) and
    therefore cannot be a standard monomial.  Since there must be $r$ standard monomials for
    $\mathbf I(V)\cap k[EV_i]$, these must by default be the monomials returned in Step~(6) of SM-A.
\end{proof}

\begin{rem}
\label{SMrem}
    At the end of Step 2 of EssGB, the set $SM_n$ is
    indeed the set of standard monomials for $\mathbf I(V)\cap k[EV_n]$
    (Corollary 2 in \cite{JS}). Furthermore, it is the set of standard
    monomials for $\mathbf I(V)$ with respect to $\prec$ (Theorem 7 in~\cite{JS}).
\end{rem}

\subsection{LT-A}

 Let $LT$ be the output returned by LT-A
and let $EV_n$ and $SM_n$ be the sets of essential variables and the
standard monomials $SM_n$ as computed in Step 2 of EssGB. Define $B$
to be the set $B=\{x_i:x_i\notin EV_n\} \cup \{x^a:supp(x^a)\subset
EV_n, x^a \notin SM_n, x^a \text{ minimal}\}$, where \emph{minimal}
means no monomial in $LT$ divides $x^a$.

\begin{lem}
\label{LTlemma}
   Let $G$ be a \gb for $\mathbf I(V)$.
    Then the leading term ideal $LT(G)$ is generated by $B$.
\end{lem}

\begin{proof}
    Since the sets of leading terms and of standard monomials for an
    ideal are mutually exclusive, by definition $B\subset LT(G)$. Let $x^a\in LT(G)$.
    If $supp(x^a)\not \subset EV_n$, then there is
    $x_i\in B$ that divides $x^a$.  Now suppose $supp(x^a)\subset EV_n$.
    Clearly $x^a\not \in SM_n$.  Since there are a finite number of divisors of $x^a$,
    there is $x^b\in B$ that divides~$x^a$. Hence, $B$ generates $LT(G)$.
\end{proof}

Note that $B$ represents the minimum set (with respect to inclusion)
of generators for $LT(G)$. In particular, no monomial $x^a\in B$
divides any other monomial in $B$. Furthermore, the set $LT$
returned by LT-A is the set $B$.

\subsection{SP-A}

We know that separators exist (see Corollary 2.14 in
\cite{robbiano}).  We also know that separators have a canonical
form.

\begin{lem}
\label{SEPlemma}
    Let $P\subset k^n$ be a set of points, $\prec$ a term order, and
    $G$ a \gb of $\mathbf{I}(P)$ with respect to $\prec$.
    The reduced separators of the points in $P$
    can be written uniquely in terms of the standard monomials in
    $SM(G)$.
\end{lem}

\begin{proof}
    Let $f$ be a separator of a point in $P$.  Since there is $p\in P$ such that
    $f(p)=1$, then $f\not \in \mathbf{I}(P)$. Hence $f$
    is a nonzero element of $R/\mathbf{I}(P)$.  As $R/\mathbf{I}(P)$
    is generated (as a $k$-vector space) by $SM(G)$, then $f$ has a
    unique $k$-linear representation in terms of the standard
    monomials which is reduced with respect to $G$.
\end{proof}

\subsection{EssGB}

\begin{lem}
\label{PLUlemma}
    For all $i$, the system $P_{i-1}L_{i-1}U_{i-1}\cdot c=b_i$ obtained during Step 2(b)
    of the execution of the algorithm EssGB has at most one solution.
\end{lem}

\begin{proof}
    Recall that $A_{i-1}:=P_{i-1}L_{i-1}U_{i-1}$ is the $(m\times
    r)$-matrix $(x^{a_j}(p_t))$ where the monomials $x^{a_1},\ldots ,x^{a_r}\in SM_{i-1}$ and $p_t$
    is the point in row $t$ of $M$.  Since the monomials $x^{a_1},\ldots
    ,x^{a_r}$ are chosen to be linearly independent, the rank
    of $A_{i-1}$ is~$r$. Hence $A_{i-1}$ has a trivial null space.
\end{proof}

Recall that the LT-A algorithm returns the minimum set $LT =B $
of generators for the leading term ideal of~$\mathbf I(V)$.

\begin{lem}
\label{GBlemma}
    A finite set $G\subset \mathbf
I(V)$ is the reduced \gb of $\mathbf I(V)$ with respect to $\prec$
if and only if
    \begin{enumerate}
        \item $G$ is monic.
        \item $\{LT(g):g\in G\}=B$ and
        \item $Supp(tail(g)) \subset SM_n$ for every $g\in
            G$.
    \end{enumerate}
\end{lem}

\begin{proof}
    Let $I=\mathbf I(V)$. If $G$ is the reduced \gb for $I$, then (1)
    holds by definition. By
    Lemma~\ref{LTlemma}, $B \subseteq \{LT(g):g\in G\}$.
    On the other hand, we cannot have different $g,h \in G$
    with the leading term of $g$ dividing the leading term
    of $h$.
    Therefore $\{LT(g):g\in G\}$ must be equal to the minimum
    set $B$ of its generators, and~(2) holds.  Moreover, if
    $x^a \notin SM_n$, then there must be some $x^b \in B$ that divides
    $x^a$, and hence $x^a$ cannot be in $Supp(tail(g))$ for any
    $g \in G$, which is equivalent to condition~(3).

    Now let $G\subset I$ be a finite set that satisfies (1)--(3).
    Let $H$ be any Gr\"obner basis, and let $f\in \mathbf I(V)$.
    Then the leading monomial of some $h \in H$ divides $LT(f)$, and
    by Lemma~\ref{LTlemma}, some $x^a \in B$ divides $LT(f)$.  Now~(2) implies that $LT(g)$ divides $LT(f)$ for some
    $g \in G$.  Thus $G$ is a Gr\"obner basis and is monic by~(1).

    Finally, let $g, h$ be different elements of $G$.  Then $LT(g)$ does not divide $LT(h)$ by minimality of $B$.
    Moreover, $LT(g)$ cannot divide any monomial in $Supp(tail((h))$,
     since by~(3) the latter monomials
    are standard monomials, while $LT(g)$ is not in $SM_n$.
\end{proof}

\begin{thm}
    Let $(G, SM_n, S)$ be the output returned by the EssGB algorithm,
    given a variety $V$ and a term order $\prec$.
    Then $G$ is the reduced \gb of $\mathbf{I}(V)$ with respect to $\prec$,
    $SM_n$
    is the set of standard
    monomials associated to $G$, and $S$ is the set of reduced separators
    of the points in $V$.
\end{thm}

\begin{proof}
    This follows from Remark \ref{SMrem} and Lemmas \ref{SEPlemma} and \ref{GBlemma}.
\end{proof}

The algorithm EssGB can be simplified for lexicographical orders.
Specifically once a monomial has been identified as a standard
monomial in Step~$i$ of EssGB, then it continues to be a standard
monomial in subsequent iterations. This property can be used to
simplify the algorithm SM-A for the case of lexicographical orders.
However, the simplification would not reduce the order of magnitude
of our worst-case run-time estimate, and we did not implement it.

\section{Complexity of the Algorithms}

\subsection{Complexity of SM-A}

Let $p=|k|$.  There are $O(pm)$ candidate monomials, which require
$O(pm^2\log(pm))$ steps to sort, assuming that comparison of two
exponents is an operation of cost $O(m)$.  The matrix $A_i$ has
$O(m\cdot pm)$ entries. Computing the row-echelon form of $A$ has
time complexity $O(pm\cdot m^2)$. Identification of the columns with
pivots is an $O(m\cdot pm)$ operation. Hence the worst-case
complexity of SM-A is
$$O(pm+pm^2\log(pm)+pm^3+pm^2)=O(pm^2\log(pm)+pm^3).$$

\subsection{Complexity of LT-A}

As there are at most $m^2$ candidate monomials, initialization of
the matrix $C$ requires $O(m^2)$ operations. The FOR loop in Step~2
is executed $O(m)$ times, similarly for the FOR loop in Step~2(a).
Checking for membership of $SM_n$ in the IF clause of Step 2(a)(i)
requires $O(m)$ operations. Checking for divisibility in
Step~2(a)(iii) can be implemented by using a look-up table, and can
be presumed to have a constant cost in each iteration of~2(a)(iii),
while creating the look-up table requires a one-time cost of
$O(m^3)$.  In all, the cost associated to Step~2 is $O(m^3)$.
 Step~3 requires $O(m^2)$ computations, while Step~4
requires $O(m^2\log m)$ computations. The last step requires $O(n)$
computations as there are at most $n-m$ inessential variables.
Overall the complexity of LT-A is
$$O(m^2+ m^3+ m^2+m^2\log m+n)=O(n+m^3).$$

\subsection{Complexity of SP-A}

Initialization of the set $S$ is a constant operation. For the FOR
loop, since we are using the PLU decomposition of the matrix $A$,
solving each of the $m$ systems in 2(a) requires $O(m^2)$ steps for
forward and backward substitution.  Maintenance of the set $S$ in
2(b) requires $m$ scalar multiplications.  Hence, the complexity of
the SP-A algorithm is $O(m)O(m^2+m)=O(m^3)$.

\subsection{Complexity of EssGB}

Initialization has cost $O(m)$. In the main FOR loop (Step~2),
the IF statement assumes that we have a linear system in PLU
form and so requires $O(m^2)$ operations for solving the system
using forward and backward substitutions. Given no solution
(entering the THEN clause), to compute the new set of standard
monomials is $O(pm^2\log(pm)+pm^3)$. Construction of the matrix
$A_i$ requires $O(m^2)$ operations since the numbers of its
rows and columns are both bounded above by $m$ and another
$O(m^3)$ to compute its PLU decomposition.  Since there are at
most $m$ essential variables, the THEN clause will only be
executed $O(m)$ times, resulting in
$$O(n)O(m^2)+O(m)O(m^2+pm^2\log(pm)+pm^3+m^2+m^3)=O(nm^2+pm^3\log(pm)+pm^4)$$
as the total cost for Step 2.

Executing Step 3 is $O(n+m^3)$, as derived above. Construction of
$b_j$ in Step 4(a) requires $O(m)$ operations, while solving the
system in 4(b) requires $O(m^2)$ operations each for forward and
backward substitution. Appending to the list $GB$ is an $O(m)$
operation. Since there are at most $n+m^2$ leading terms, the total
cost of Step 4 is $O(n+m^2)O(m^2)=O(nm^2+m^4)$.

Executing Step 5 is $O(m^3)$, as derived above. Thus the worst-case
complexity of the EssGB algorithm is
$$O(m)+O(nm^2+pm^3\log(pm)+pm^4)+O(n+m^3)+O(nm^2+m^4)+O(m^3)$$
$$=O\left(pm^3\log(pm)+pm^4+nm^2\right).$$
If we assume $p$ to be fixed, then the complexity can be reduced to
$O(nm^2+m^4)$.  For the applications to biological data where $n\gg
m$, the complexity is dominated by $O(nm^2)$.

\section{Performance of the EssGB Algorithm}

We compared the run-times of the algorithms EssGB, EssBM, and
BM on randomly generated varieties containing~$m$ points in
$k^n$, where $k$ is a finite field of the form $\mathbb
Z/p\mathbb Z$. We performed this comparison in \emph{Macaulay
2}, version 0.9.97, where each algorithm has been implemented.

We generated $r=10$ affine varieties for changing values of $p$,
$n$, and $m$.  Since the algorithms require specification of a term
order, we consider this to be parameter as well.  The table below
lists the values we used for this comparative study.

\bigskip
\noindent%
\begin{center}
\begin{tabular}{|l|l|}
  \hline
  Parameters & Values \\ \hline
  $p$  = cardinality of $k$ & \{5, 101\} \\
  $n$  = number of variables & \{100, 200, 300\} \\
  $m$  = number of points & \{5, 10, 15\} \\
  $\prec$  = term order & \{Lex, GRevLex\} with default variable order  \\
  \hline
\end{tabular}
\end{center}
\bigskip

For $1 \leq i \leq r$, the $i$-th variety consists of $nr(i)$
randomly generated points, where $$nr(i) = \frac{m}{5} \left\lceil
\frac{r-i+1}{2}\right\rceil.$$  The remaining $m-nr(i)$ points were
generated using random homogenous linear polynomials
$g_1,\ldots,g_{m-nr(i)}$, where $g_j \in k[y_1,\ldots,y_{nr(i)}]$.
To generate the $j$-th new point $p_j$,  the coordinates of $p_j$
are computed individually; that is, for $1\leq \ell \leq n$
$$p_{j\ell} := g_j\left(p_{1\ell},\ldots,p_{nr(i)\ell}\right).$$

Note that for $i=1,2$, all points are randomly generated. This will
result, with probability very close to one, in a variety where the
points are in \emph{general position}, that is, there are no linear
dependencies among the points (\cite{harris}, pg. 7).  In the runs
for $i = 3, \dots, 10$, the enforced randomly chosen linear
dependencies ensure that the linear span of the generated variety
will have dimension $\leq nr(i)$ (and equal to $nr(i)$ with
probability close to one). This choice of test data allowed us to
compare run-times of the three algorithms on ideals of varieties
with different geometric properties.

\begin{figure}[h]
\centering
  \includegraphics[width=5in]{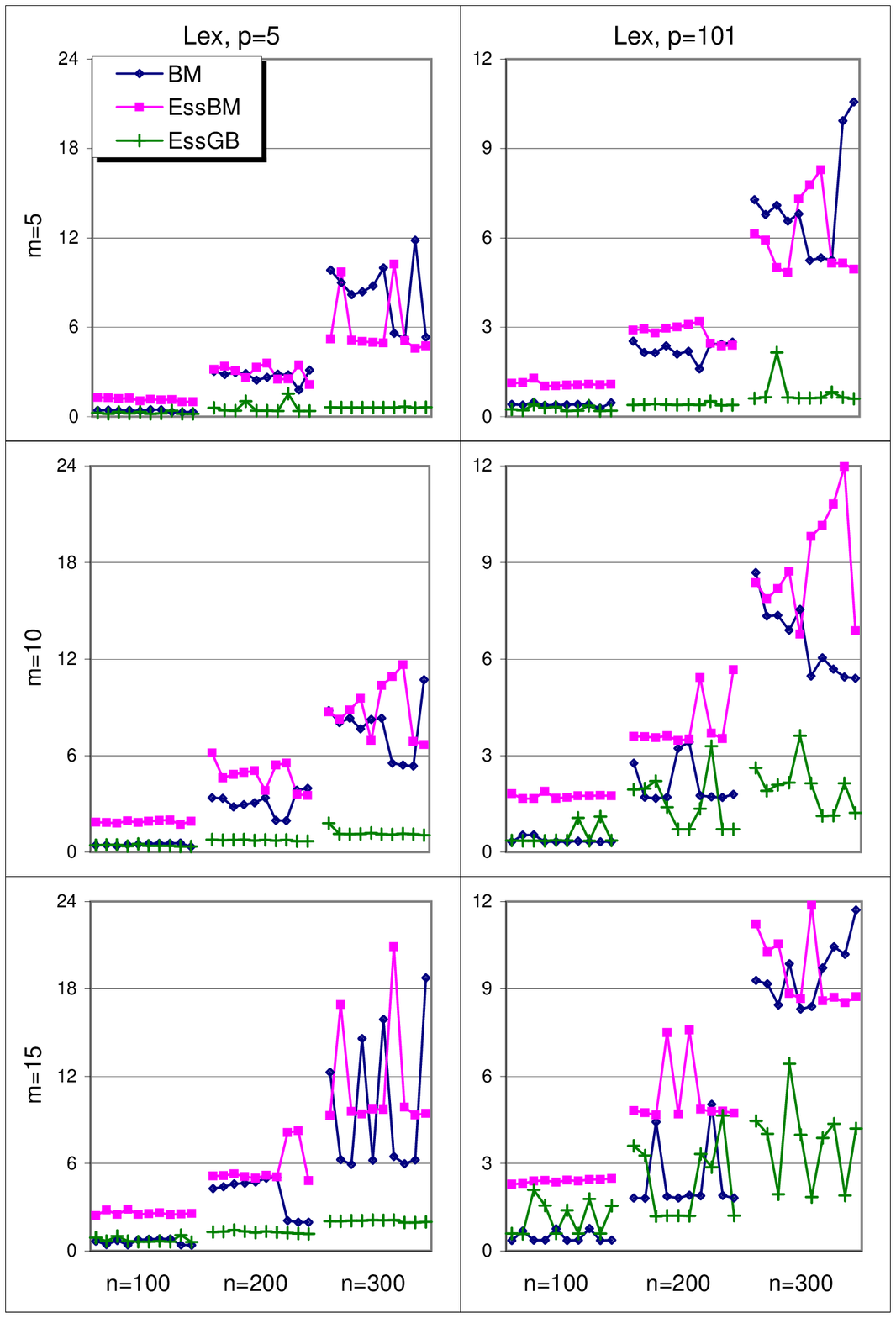}\\
  \caption{Run-times for the algorithms BM, EssBM, and EssGB for $p \in \{5,
101\}$, $m \in \{5, 10, 15\}$, and $n\in \{100,200,300\}$ with a
default Lex order. } \label{lex-plots}
\end{figure}

\begin{figure}[h]
\centering
  \includegraphics[width=5in]{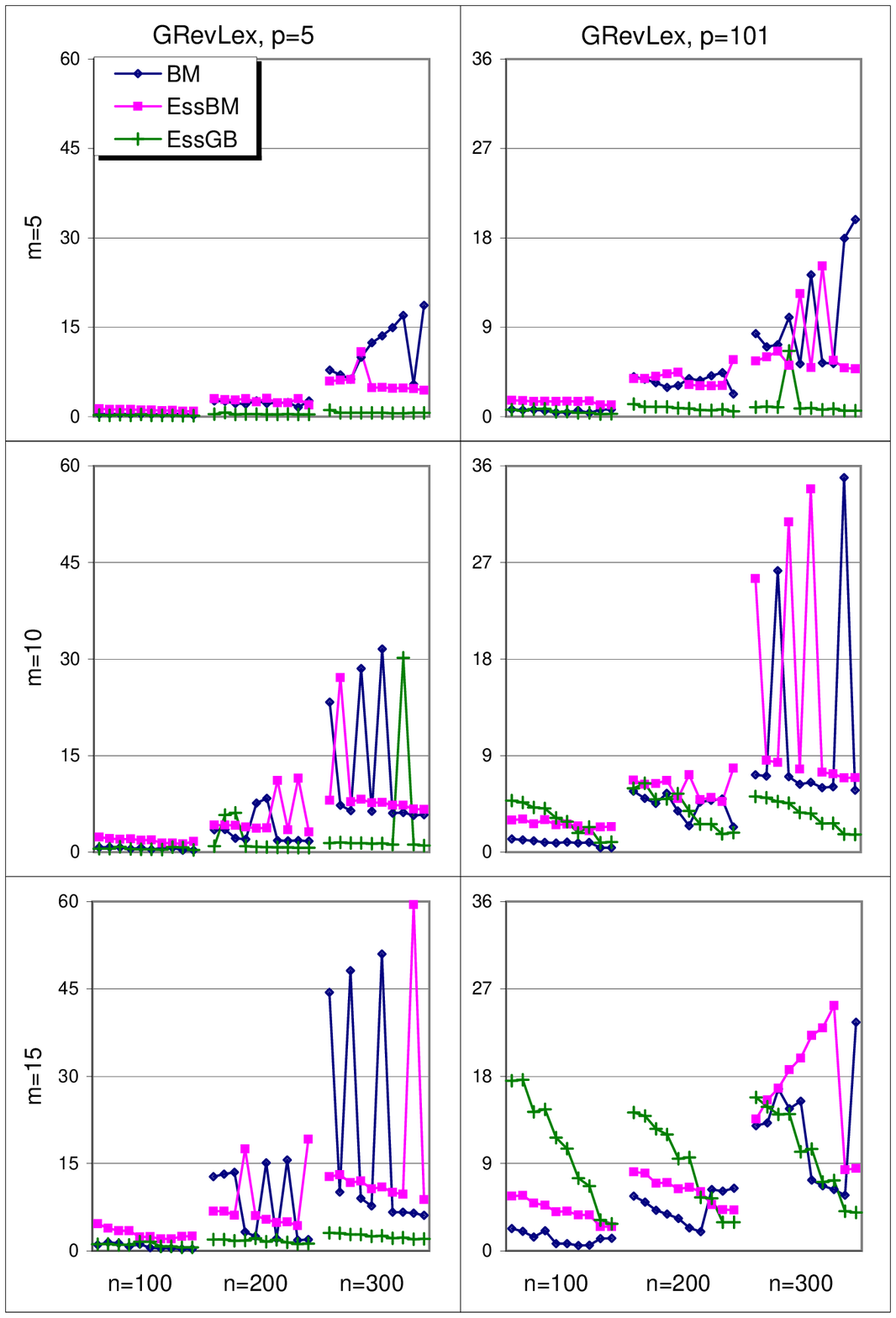}\\
  \caption{Run-times for the algorithms BM, EssBM, and EssGB for $p \in \{5,
101\}$, $m \in \{5, 10, 15\}$, and $n\in \{100,200,300\}$ with a
default GRevLex order. } \label{grevlex-plots}
\end{figure}

We applied the three algorithms to each of the generated varieties.
The run-time results are displayed in Figures \ref{lex-plots} and
\ref{grevlex-plots}.

\section{Discussion}

Recently, \gbs have been used as a selection tool in applications to
molecular biology \cite{LS, DJLS}.
 In these applications,
the number of data points~$m$ tends to be significantly smaller than
the number of variables~$n$. The computation of \gbs constitutes a
bottleneck for overall feasibility of these calculations.  The
primary motivation for our paper was to find an algorithm that
optimizes run-time in the case when $m \ll n$.

 The time complexity of the standard BM
algorithm has been reported in the literature as quadratic in the
number of variables $n$ and cubic in the number of points~$m$
\cite{BM}.  This makes it too slow for the applications mentioned in
the preceding paragraph.  In \cite{JS}, we developed an algorithm
EssBM that has a provable worst-case time complexity of $O(nm^3 +
m^6)$ for a fixed finite field~$k$. For the algorithm EssGB
presented here, we can improve this worst-case estimate to $O(nm^2 +
m^4)$ for a fixed finite field~$k$.  The reduction from quadratic to
linear scaling in the run-time was achieved in both EssBM and EssGB
by first identifying the set of essential variables in a single loop
of length~$n$, and performing the most expensive steps of the
computation only for these essential variables.  While EssBM still
uses BM as a subroutine on the reduced set of variables, EssGB
eliminates calls to BM altogether and computes all relevant objects
by solving systems of $k$-linear equations.  The coefficient
matrices used in these equations change only when a new essential
variable is encountered.  This allows us to use PLU decompositions
to reduce the cost to $O(m^2)$ in all but~$m$ of the~$n$ steps of
the main loop, and our overall worst-case estimate follows.

Based on this estimate, one would expect our algorithm to be
significantly faster than both BM and EssGB when $m \ll n$.  We
tested this prediction for randomly generated varieties, with $|V|=
m \in \{5, 10, 15\}$.  We tested the algorithm  on
varieties that were generated totally randomly, which should ensure
that the points will almost certainly be in general position, and on
varieties where an increasing number of the points were
expressed as  linear
combinations of previously defined random points. In order to ensure
that we have enough different linear combinations of two points, the
smallest field for which we tested our algorithm was $\mathbb
Z/5\mathbb Z$.  We also run tests for the rather large field
$\mathbb Z/{101}\mathbb Z$.

Our test runs neatly confirm that our algorithm EssGB has comparable
performance with BM when $n = 100$, and significantly outperforms
the latter when $n = 300$.  The single exception are the simulations
where $m = 15$, $p = 101$, and a GrevLex term order $\prec$ was
used. In these simulations the performance of our algorithm becomes
only comparable to that of BM when $n = 300$.  However, the general
pattern still holds: The more variables, the better EssGB performs
relative to BM.

We also observed that in general the run-times of EssGB are more
consistent for different varieties under the same parameter settings
than those for BM or EssBM.  The only exception here are the
experiments with  $m = 10, 15$, $p = 101$, and a Lex term order
$\prec$, where similar magnitudes of run-time fluctuations were
observed for all three algorithms. The experiments with GrevLex term
orders and $p= 101$ also show a significant decrease of the run-time
of EssGB when the number of linear dependencies among the points in
the variety increases.

Our simulations do not in general show an advantage of our previous
algorithm EssBM over BM, although EssBM clearly does become more
competitive with BM as the number~$n$ of variables increases.
Previous experiments reported in \cite{JS} had shown that EssBM
outperforms BM when the number of variables starts exceeding 200.
However, these experiments were run in implementation 0.9.8 of
\emph{Macaulay~2,} while the simulations presented here were run on
version 0.9.97.  We noticed a significant speedup of the run-times
for both BM and EssBM between both versions; it was relatively
larger for BM.

In summary, both our theoretical run-time estimates and the test
runs reported here indicate that EssGB would be the algorithm of
choice if Gr\"obner bases are to be found for a variety $V$ in
$k[x_1, \ldots , x_n]$ such that $|V| = m \ll n$.

\section{Acknowledgements}

The authors would like to thank Xiaoping A. Shen and Luis Garc\'ia
for their helpful suggestions and discussions. This material is
based upon work supported by the National Science Foundation under
Agreement No. 0112050.

\bibliographystyle{amsplain}
\bibliography{mod-bm}

\providecommand{\bysame}{\leavevmode\hbox to3em{\hrulefill}\thinspace}
\providecommand{\MR}{\relax\ifhmode\unskip\space\fi MR }
\providecommand{\MRhref}[2]{%
  \href{http://www.ams.org/mathscinet-getitem?mr=#1}{#2}
}
\providecommand{\href}[2]{#2}
\begin{thebibliography}{10}

\bibitem{abbott}
J.~Abbott, A.~Bigatti, M.~Kreuzer, and L.~Robbiano, \emph{Computing ideals of
  points}, Journal of Symbolic Computation \textbf{30} (2000), no.~4, 341--356.

\bibitem{cerlienco}
L.~Cerlienco and M.~Mureddu, \emph{From algebraic sets to monomial linear bases
  by means of combinatorial algorithms}, Discrete Mathematics \textbf{139}
  (1995), 73–--87.

\bibitem{DJLS}
E.~Dimitrova, A.~Jarrah, R.~Laubenbacher, and B.~Stigler, \emph{A {G}r\"obner
  fan method for biochemical network modeling}, ISSAC '07: Proceedings of the
  2007 international symposium on Symbolic and algebraic computation (New York,
  NY, USA), ACM, 2007, pp.~122--126.

\bibitem{farr}
J.~Farr and S.~Gao, \emph{Computing {G}r{\"o}bner bases for vanishing ideals of
  finite sets of points}, Applied Algebra, Algebraic Algorithms and
  Error-Correcting Codes: 16th International Symposium, AAECC-16 (M.~Fossorier,
  H.~Imai, S.~Lin, and A.~Poli, eds.), Lecture Notes in Computer Science, vol.
  3857, Springer Berlin, 2006, pp.~118--127.

\bibitem{harris}
J.~Harris, \emph{Algebraic geometry: A first course}, 1st ed., Graduate Texts
  in Mathematics, vol. 133, Springer-Verlag, New York, 1992.

\bibitem{JS}
Winfried Just and Brandilyn Stigler, \emph{Computing {G}r\"obner bases of
  ideals of few points in high dimensions}, Communications in Computer Algebra
  \textbf{40} (2006), no.~3, 65--96.

\bibitem{LS}
R.~Laubenbacher and B.~Stigler, \emph{A computational algebra approach to the
  reverse engineering of gene regulatory networks}, Journal of Theoretical
  Biology \textbf{229} (2004), 523--537.

\bibitem{lederer}
M.~Lederer, \emph{The vanishing ideal of a finite set of closed points in
  affine space}, Available at {\tt http://arxiv.org/abs/math/0604133}, 2006.

\bibitem{marinari}
M.~Marinari, H.~M. M\"oller, and T.~Mora, \emph{Gr\"obner bases of ideals
  defined by functionals with an application to ideals of projective points},
  Applicable Algebra in Engineering, Communication and Computing \textbf{4}
  (1993), 103--145.

\bibitem{BM}
H.~M. M\"oller and B.~Buchberger, \emph{The construction of multivariate
  polynomials with preassigned zeroes}, Computer Algebra: EUROCAM '82
  (J.~Calmet, ed.), Lecture Notes in Computer Science, vol. 144, Springer
  Berlin, 1982, pp.~24--31.

\bibitem{mora}
T.~Mora and L.~Robbiano, \emph{Points in affine and projective spaces},
  Computational Algebraic Geometry and Commutative Algebra, Cortona-91
  (D.~Eisenbud and L.~Robbiano, eds.), Symposia Mathematica, vol.~34, Cambridge
  University Press, 1993, pp.~106--150.

\bibitem{robbiano}
L.~Robbiano, \emph{Gr\"{o}bner bases and statistics}, Gr\"{o}bner Bases and
  Applications (New York) (B.~Buchberger and F.~Winkler, eds.), London
  Mathematical Society Lecture Notes Series, vol. 251, Cambridge University
  Press, 1998, pp.~179--204.

\bibitem{yeung}
M.~K.~S. Yeung, J.~Tegn\'er, and J.~Collins, \emph{Reverse engineering gene
  networks using singular value decomposition and robust regression}, PNAS
  \textbf{99} (2002), no.~9, 6163--6168.

\end{thebibliography}

\end{document}